\documentclass{article}
\usepackage{amsmath,amssymb,amsthm,graphicx,epsfig,subfigure,float,url}
\usepackage[colorlinks=true]{hyperref}
\usepackage{pdfsync}
\usepackage{esint}

\usepackage{verbatim}
\usepackage[applemac]{inputenc}
 
\usepackage[usenames,dvipsnames]{pstricks}
\usepackage{pst-grad} 
\usepackage{pst-plot} 

\def\R{\textrm{I\kern-0.21emR}}
\def\N{\textrm{I\kern-0.21emN}}
\def\1{1\textnormal{\kern -0.21emI}}

\renewcommand{\geq}{\geqslant}
\renewcommand{\leq}{\leqslant}

\newtheorem{theorem}{Theorem}

\theoremstyle{definition}
\theoremstyle{definition}

\title{Optimal control of resources for species survival}
\author{Idriss Mazari\footnote{Sorbonne Universit\'es, UPMC Univ Paris 06, UMR 7598, Laboratoire Jacques-Louis Lions, F-75005, Paris, France (\texttt{idriss.mazari@upmc.fr}).}
\and Gr\'egoire Nadin\footnote{CNRS, Universit\'e Pierre et Marie Curie (Univ. Paris 6), UMR 7598, Laboratoire Jacques-Louis Lions, F-75005, Paris, France ({\tt gregoire.nadin@upmc.fr}).}
	\and Yannick Privat\footnote{IRMA, Universit\'e de Strasbourg, CNRS UMR 7501, 7 rue Ren\'e Descartes, 67084 Strasbourg, France ({\tt yannick.privat@unistra.fr}).}
}

\date{}

\begin{document}

\maketitle

\begin{abstract}
Consider a species whose population density solves the steady diffusive logistic equation in a heterogeneous environment modeled with the help of a spatially non constant coefficient standing for a resources distribution in a given box. We look at maximizing the total population size with respect to resources distribution, under some biologically relevant constraints. Assuming that the diffusion rate of the species is large enough, we prove that any optimal configuration is the characteristic function of a domain standing for the resources location. Moreover, we highlight that optimal configurations look {\it concentrated} whenever the diffusion rate is large enough. 
In the one-dimensional case, this problem is deeply analyzed, and for large diffusion rates, all optimal configurations are exhibited. 
\end{abstract}

\section{Setting of the problem}
Let ${\Omega}=(0,1)^n$. We  investigate an optimal control problem arising in population dynamics. 
Consider the population density ${\theta_{m,\mu}}$ of a given species evolving in ${\Omega}$, assumed to solve the so-called steady logistic-diffusive equation writing 
\begin{equation}\tag{LDE}\label{LDE} \left\{
\begin{array}{ll}
\mu {\Delta} {\theta_{m,\mu}}(x)+(m(x)-{\theta_{m,\mu}}(x)){\theta_{m,\mu}}(x) =0 & x\in {\Omega},\\
\frac{\partial {\theta_{m,\mu}}}{\partial \nu} =0 &x\in \partial {\Omega},
\end{array}
\right.\end{equation} 
where $m$ is a bounded function of $\Omega$ standing for the resources distribution and $\mu>0$ stands for the species velocity also called {\it diffusion rate}. From a biological point of view, the real number $m(x)$ is a measure of the resources available at $x$ of the habitat ${\Omega}$ and can be seen as the local intrinsic growth rate of species at location $x\in \Omega$.

The pioneering works by Fisher \cite{Fisher}, Kolmogorov-Petrovski-Piskounov \cite{KPP} and Skellam \cite{SkellamRandom} sparked a new interest in the study of the heterogeneity influence of the environment on the growth of a population densities. In the present work, we look at investigating the influence of heterogeneous environment on the density population. To this aim, for a given $\mu>0$, we address the optimal control problem of {\it maximizing the total population size}, given by 
\begin{equation}\label{def:Fmu}
F_\mu(m)= \int_{\Omega} {\theta_{m,\mu}},
\end{equation}
where $ {\theta_{m,\mu}}$ denotes the solution of equation \eqref{LDE}. 
In the framework of population dynamics, the density $\theta_{m,\mu}$ solving Equation \eqref{LDE} can be interpreted as a steady state associated to the following evolution equation
\begin{equation}\tag{LDEE}\label{LDEE} \left\{
\begin{array}{ll}
 \frac{\partial u}{\partial t}(t,x) = \mu {\Delta} u(t,x)+u(t,x)(m(x)-u(t,x)) & t>0, \ x\in {\Omega}\\
\frac{\partial u}{\partial \nu}(t,x)=0& t>0, \ x\in \partial {\Omega}\\ 
u(0,x)=u^0(x) & x\in {\Omega} \\
\end{array}
\right.\end{equation}
modeling the spatiotemporal behavior of a population density $u$ in a domain ${\Omega}$ with the spatially heterogeneous resource term $m$.
More precisely, provided that $\int_\Omega m>0$, the steady state ${\theta_{m,\mu}}$ is globally asymptotically stable meaning that, for any nonnegative and nonzero function $u^0\in H^1({\Omega})$, $u(t,\cdot)$ converges uniformly to ${\theta_{m,\mu}}$ as $t\to +\infty$.

From the biological point of view, it is relevant to assume that the density $m$ is nonnegative and bounded by a positive constant $\kappa$ and that $\fint_{\Omega} m$, the total amount of resources, is fixed. This leads to address the issue:
\begin{center}
{\it For given diffusivity and total amount of resources, which weight $m$ maximizes the total population size among all uniformly bounded elements ?}
\end{center}
This issue rewrites more precisely
\begin{center}
\fbox{
\begin{minipage}{15cm}
\begin{center}
{\bf Optimal design problem. }\textit{
Fix $n\in \textrm{I\kern-0.21emN}^*$, $\mu>0$, $\kappa>0$, $m_0\in (0,\kappa)$. 
We consider the optimization problem
\begin{equation}\tag{$\mathcal{P}_\mu$}\label{TSOM}
\sup_{m\in \mathcal{M}_{	\kappa,m_0}({\Omega})} F_\mu(m)\quad \text{with}\quad \mathcal M_{\kappa,m_0}({\Omega}):=\left\{m\in L^\infty({\Omega})\, , 0\leq m\leq \kappa \text{ a.e and }\fint_{\Omega} m=m_0\right\}.
\end{equation} }
\end{center}\end{minipage}}
\end{center}

It is notable that similar issues for linear models have been addressed for instance in \cite{JhaPorru,KaoLouYanagida,LamboleyLaurainNadinPrivatProperties,LouSome,MR2281509,Roques-Hamel}.

\section{Optimal resources sets}
This section is devoted to stating two important properties of solutions of \eqref{TSOM}, whose proof can be found in \cite{MNP}.

 \paragraph{Property no. 1: pointwise constraints, bang-bang property.}
We investigate hereafter the {\it bang-bang} character of maximizers for Problem \eqref{TSOM},  in other words, whether a solution $m^*$ is equal to 0 or $\kappa$ a.e. in $\Omega$. Such an issue is of importance, since dealing with {\it bang-bang} functions enables the use of adapted numerical approaches. 

\begin{theorem}[\cite{MNP}, Theorem 1]\label{TheoremePrincipal}
Let $n\in \textrm{I\kern-0.21emN}^*$, $\mu>0$, $\kappa>0$, $m_0\in (0,\kappa)$. There exists a positive number $\mu^*=\mu^*({\Omega},\kappa,m_0)$ such that, for every $\mu\geq \mu^*$ every  maximizer of Problem \eqref{TSOM} is  bang-bang. As a consequence, $m^*=\kappa \chi_E$, where $\chi_E$ is the characteristic function of a (measurable) resources set $E$. 
\end{theorem}

The main difficulty for establishing this result rests upon the facts that $\theta_{m,\mu}$ solves a nonlinear PDE and the criterion $F_\mu$ does not derive from an energy, making the exploitation of adjoint approaches intricate. 

 \paragraph{Property no. 2: concentration-fragmentation of maximizers.} It is well-known (see e.g. \cite{BHR,LamboleyLaurainNadinPrivatProperties}) that concentrating resources, meaning that the resources distribution $m$ is decreasing in each direction, favors the survival of species. On the contrary, we will say that a resources set is fragmented whenever it is disconnected. On Figures \ref{fig:3}-\ref{fig:4} below, ${\Omega}$ is a square, and the intuitive notion of {\it concentration-fragmentation} of resources distributions is illustrated.
\begin{center}
\begin{figure}[h!]
\centering
\begin{minipage}{6cm}
\includegraphics[height=50mm]{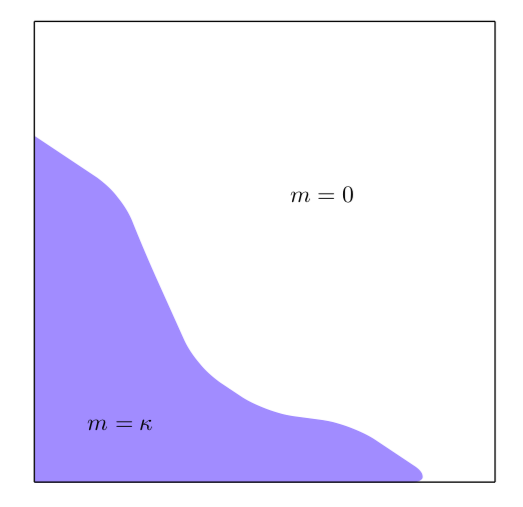}
\caption{Concentrated distribution.}
\label{fig:3}
\end{minipage}
\hspace{1cm}
\begin{minipage}{6cm}
\includegraphics[height=50mm]{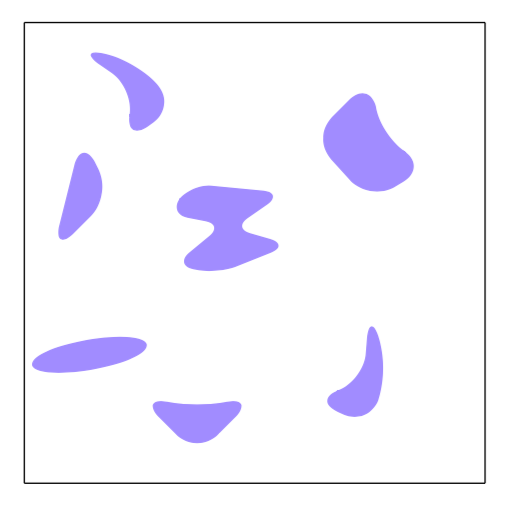}
\caption{Fragmented distribution.}
\label{fig:4}
\end{minipage}\ \\
\begin{minipage}{9cm}
\includegraphics[width=5cm]{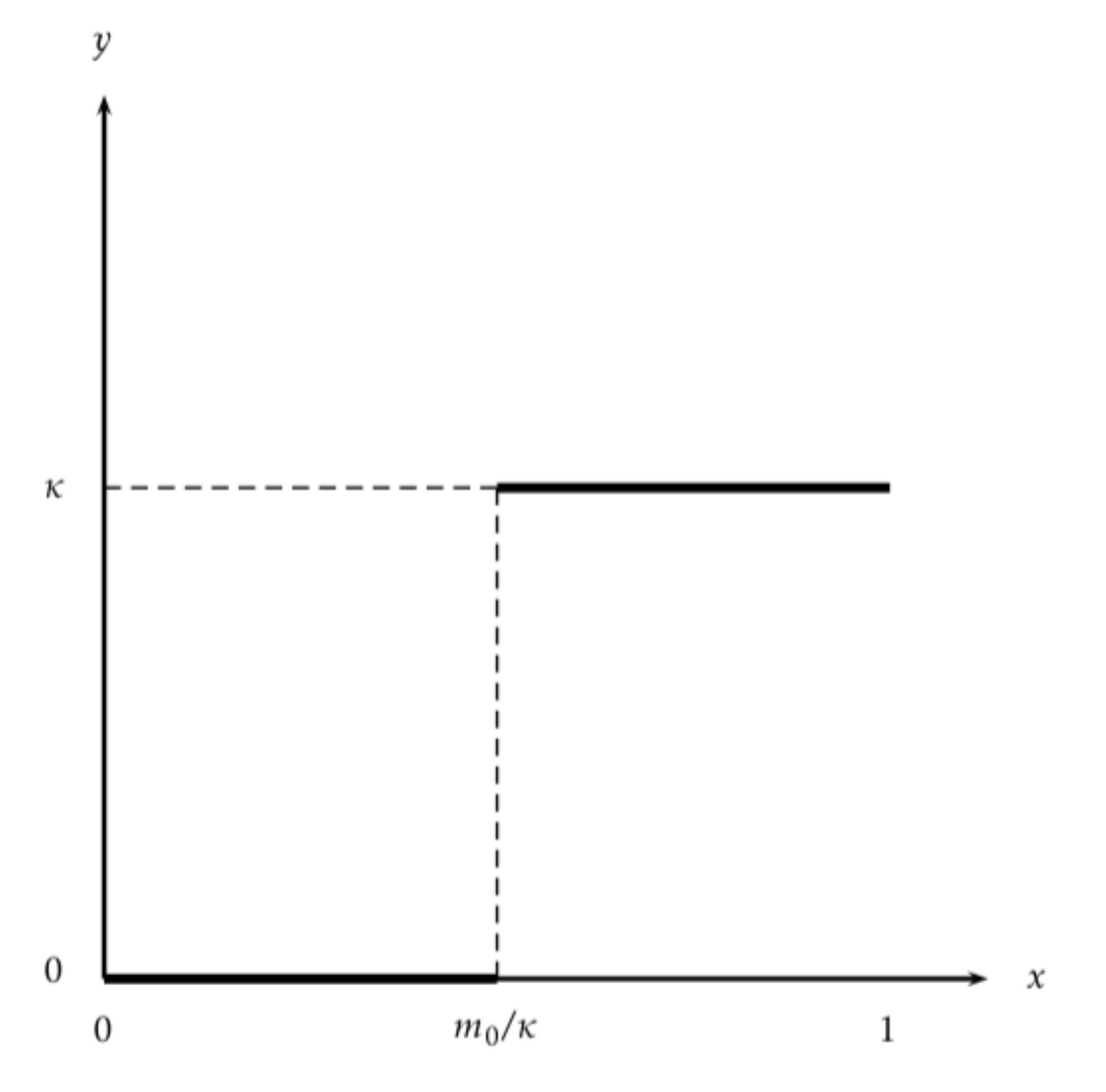}
\caption{A solution of \eqref{TSOM} in 1D, for large diffusivities $\mu$.}
\label{fig:5}
\end{minipage}
\end{figure}
\end{center}
\vspace{-0.5cm}
\begin{theorem}[\cite{MNP}, Theorem 2]\label{Theoreme2}
Any family of maximizers $\{m_\mu\}_{\mu >0}$ converges in $L^1(\Omega)$ to the characteristic function of a set $E$ which is concentrated (meaning that its characteristic function $\chi_E$ is monotone with respect to each space variable).
\end{theorem}

In the one-dimensional case, we also prove that {if the diffusivity is large enough, there are only two maximizers}, that are simple crenels meeting either the left or the right boundary (see Fig. \ref{fig:5} above).
Finally, in the one-dimensional case, we obtain a surprising result: \textbf{fragmentation may be better than concentration for small diffusivities}.


\paragraph{Acknowledgement.}
The authors were partially supported by the Project ``Analysis and simulation of optimal shapes - application to lifesciences'' of the Paris City Hall.

\end{document}